\documentclass[10pt]{article}

\usepackage{amscd,amsmath, amssymb, fancyhdr, url, graphicx}

\numberwithin{equation}{section}

\newcommand{\version}{version 2.1,\ \ Oct 26, 2014}

\makeatletter
\def\x@arrow{\DOTSB\Relbar}
\def\xlongrightarrowfill@{\arrowfill@\relbar\relbar\longrightarrow}
\newcommand{\xlongrightarrow}[2][]{%
        \ext@arrow 0099\xlongrightarrowfill@{#1}{#2}}
\makeatother

\def\eqref#1{(\ref{#1})}

\newcommand{\arrow}{{\:\longrightarrow\:}}
\newcommand{\Z}{{\Bbb Z}}
\def\C{{\Bbb C}}
\def\P{{\Bbb P}}

\newcommand{\R}{{\Bbb R}}
\newcommand{\Q}{{\Bbb Q}}
\renewcommand{\H}{{\Bbb H}}

\def\1{\sqrt{-1}\:}

\newcommand{\cntrct}                % contraction with a vector field
{\hspace{2pt}\raisebox{1pt}{\text{$\lrcorner$}}\hspace{2pt}}

\renewcommand{\tilde}{\widetilde}
\renewcommand{\bar}{\overline}
\renewcommand{\phi}{\varphi}
\renewcommand{\epsilon}{\varepsilon}
\renewcommand{\geq}{\geqslant}

\newcommand{\Teich}{\operatorname{\sf Teich}}

\newcommand{\Comp}{\operatorname{\sf Comp}}
\newcommand{\NS}{\operatorname{\sf NS}}
\newcommand{\Per}{\operatorname{\sf Per}}
\newcommand{\Perspace}{\operatorname{{\Bbb P}\sf er}}

\newcommand{\Kah}{\operatorname{Kah}}
\newcommand{\Amp}{\operatorname{Amp}}

\newcommand{\Bir}{\operatorname{Bir}}

\newcommand{\Pos}{\operatorname{Pos}}
\newcommand{\Hol}{\operatorname{Hol}}
\newcommand{\Aut}{\operatorname{Aut}}

\newcommand{\Hdg}{{\operatorname{\sf Hdg}}}

\newcommand{\Diff}{\operatorname{\sf Diff}}
\newcommand{\Supp}{\operatorname{Supp}}

\renewcommand{\Im}{\operatorname{Im}}

\newcounter{Mycounter}[section]
\newcounter{lemma}[section]
\renewcommand{\thelemma}{{Lemma \thesection.\arabic{lemma}}}
\newcommand{\lemma}{%
    \setcounter{lemma}{\value{Mycounter}}
    \refstepcounter{lemma}
    \stepcounter{Mycounter}
    {\noindent \bf \thelemma:\ }}

\newcounter{claim}[section]
\renewcommand{\theclaim}{{Claim \thesection.\arabic{claim}}}
\newcommand{\claim}{%
    \setcounter{claim}{\value{Mycounter}}
    \refstepcounter{claim}
    \stepcounter{Mycounter}
    {\noindent \bf \theclaim:\ }}

\newcounter{sublemma}[section]
\newcounter{corollary}[section]
\renewcommand{\thecorollary}{{Corollary \thesection.\arabic{corollary}}}
\newcommand{\corollary}{%
    \setcounter{corollary}{\value{Mycounter}}
    \refstepcounter{corollary}
    \stepcounter{Mycounter}
    {\noindent \bf \thecorollary:\ }}

\newcounter{theorem}[section]
\renewcommand{\thetheorem}{{Theorem \thesection.\arabic{theorem}}}
\newcommand{\theorem}{%
    \setcounter{theorem}{\value{Mycounter}}
    \refstepcounter{theorem}
    \stepcounter{Mycounter}
    {\noindent \bf \thetheorem:\ }}

\newcounter{conjecture}[section]
\renewcommand{\theconjecture}{{Conjecture \thesection.\arabic{conjecture}}}
\newcommand{\conjecture}{%
    \setcounter{conjecture}{\value{Mycounter}}
    \refstepcounter{conjecture}
    \stepcounter{Mycounter}
    {\noindent \bf \theconjecture:\ }}

\newcounter{proposition}[section]
\renewcommand{\theproposition}
      {{Proposition \thesection.\arabic{proposition}}}
\newcommand{\proposition}{%
    \setcounter{proposition}{\value{Mycounter}}
    \refstepcounter{proposition}
    \stepcounter{Mycounter}
    {\noindent \bf \theproposition:\ }}

\newcounter{definition}[section]
\renewcommand{\thedefinition}
      {{Definition~\thesection.\arabic{definition}}}
\newcommand{\definition}{%
    \setcounter{definition}{\value{Mycounter}}
    \refstepcounter{definition}
    \stepcounter{Mycounter}
    {\noindent \bf \thedefinition:\ }}

\newcounter{example}[section]
\newcounter{remark}[section]
\renewcommand{\theremark}{{Remark \thesection.\arabic{remark}}}
\newcommand{\remark}{%
    \setcounter{remark}{\value{Mycounter}}
    \refstepcounter{remark}
    \stepcounter{Mycounter}
    {\noindent \bf \theremark:\ }}

\newcounter{problem}[section]
\newcounter{question}[section]
\makeatletter

\@addtoreset{equation}{section} \@addtoreset{footnote}{section}
\makeatother

\def\blacksquare{\hbox{\vrule width 5pt height 5pt depth 0pt}}
\def\endproof{\blacksquare}

\begin{document}

%%%%%%%%%%%%%%%%%%%%%%%%%%%%%%%%%%%%%%%%%%%%%%%%%%%%%%%%%%%%
\begin{center}
{\LARGE\bf
Morrison-Kawamata cone conjecture for hyperk\"ahler manifolds\\[4mm]
}
%%%%%%%%%%%%%%%%%%%%%%%%%%%%%%%%%%%%%%%%%%%%%%%%%%%%%%%%%%%%

Ekaterina Amerik\footnote{Partially supported by 
RSCF, grant number 14-21-00053.},  
Misha
Verbitsky\footnote{Partially supported by RSCF, grant number 14-21-00053.


}

\end{center}

%%%%%%%%%%%%%%%%%%%%%%%%%%%%%%%%%%%%%%%%%%%%%%%%
{\small \hspace{0.15\linewidth}
\begin{minipage}[t]{0.7\linewidth}
{\bf Abstract} \\
Let $M$ be a simple holomorphically symplectic manifold,
that is, a simply connected compact holomorphically symplectic manifold
of K\"ahler type with $h^{2,0}=1$. Assuming $b_2(M)\neq 5$, we prove that the group
of holomorphic automorphisms of $M$ acts on the set of faces 
of its K\"ahler cone with finitely many orbits.
This statement is known as Morrison-Kawamata cone conjecture for
hyperk\"ahler manifolds. As an implication, we show that
any hyperk\"ahler manifold has only finitely many non-equivalent
birational models. The proof is based on the following observation,
proven with ergodic theory. Let $M$ be a complete Riemannian
orbifold of dimension at least three, constant negative curvature 
and finite volume, and $\{S_i\}$ an infinite set of 
complete, locally geodesic hypersurfaces.
Then the union of $S_i$ is dense in $M$.
\end{minipage}
}
%%%%%%%%%%%%%%%%%%%%%%%%%%%%%%%%%%%%%%%%%%%%%%%%

\tableofcontents

%%%%%%%%%%%%%%%%%%%%%%%%%%%%%%%%%%%%%%%%%%%%%%%%

\section{Introduction}

%%%%%%%%%%%%%%%%%%%%%%%%%%%%%%%%%%%%%%%%%%%%%%%%

%%%%%%%%%%%%%%%%%%%%%%%%%%%%%%%%%%%%%%%%%%%%%%%%%%%%%%%%%%%%
\subsection{K\"ahler cone and MBM classes}
\label{_Kahler_cone_MBM_intro_Subsection_}
%%%%%%%%%%%%%%%%%%%%%%%%%%%%%%%%%%%%%%%%%%%%%%%%%%%%%%%%%%%%

Let $M$ be a hyperk\"ahler manifold, that is, a compact,
holomorphically symplectic K\"ahler manifold. We assume that
$\pi_1(M)=0$ and $H^{2,0}(M)=\C$: the general case reduces to this by Bogomolov decomposition
(\ref{_Bogo_deco_Theorem_}). Such hyperk\"ahler manifolds
are known as simple hyperk\"ahler manifolds, or IHS (irreducible holomorphic symplectic) manifolds.
The known examples of such manifolds are
deformations of punctual Hilbert scheme of K3 surfaces, deformations of 
generalized Kummer
varieties and two sporadic ones discovered by O'Grady.
In \cite{_AV:MBM_} we gave a description of the K\"ahler cone
of $M$ in terms of a set of cohomology classes $S\subset H^2(M,\Z)$ 
called {\bf MBM classes} (\ref{mbm}). This set is of topological
nature, that is, it depends only on the deformation type of $M$.

Recall that on the second cohomology of a hyperk\"ahler manifold, there is
an integral quadratic form $q$, called the Beauville-Bogomolov-Fujiki 
form (see section \ref{prelim} for details). This form is 
of signature $(+,-, \dots, -)$ on $H^{1,1}(M)$.
Let $\Pos\subset H^{1,1}(M)$ be the positive cone, and $S(I)$ the set of
all MBM classes which are of type (1,1) on 
$M$ with its given complex structure $I$. Then the K\"ahler cone
is a connected component of $\Pos\backslash S(I)^\bot$, 
where $S(I)^\bot$ is the union of the orthogonal complements
to all $z\in S(I)$.

%%%%%%%%%%%%%%%%%%%%%%%%%%%%%%%%%%%%%%%%%%%%%%%%%%%%%%%%%%%%%%%%%%%%%%%%
\subsection{Morrison-Kawamata cone conjecture for
  hyperk\"ahler manifolds}
%%%%%%%%%%%%%%%%%%%%%%%%%%%%%%%%%%%%%%%%%%%%%%%%%%%%%%%%%%%%%%%%%%%%%%%%

The Morrison-Kawamata cone conjecture for Calabi-Yau manifolds was stated 
in \cite{_Morrison:Beyond_}.
For K3 surfaces it was already known since mid-eighties by the work of Sterk 
\cite{Sterk}.
 Kawamata in \cite{_Kawamata:cone_} proved the conjecture 
for Calabi-Yau threefolds admitting a holomorphic fibration 
over a positive-dimensional base.

In this paper, we concentrate on the following version of the cone conjecture
(see Subsection \ref{_ample_Subsection_} 
for its relation to the classical one, formulated 
for the ample cone of a projective variety).

\hfill

%%%%%%%%%%%%%%%%%%%%%%%%%%%%%%%%%%%%%%%%%%%%%%%
\definition\label{_faces_Definition_}
Let $M$ be a compact, K\"ahler manifold,
$\Kah\subset H^{1,1}(M,\R)$ the K\"ahler cone,
and $\overline\Kah$ its closure in $H^{1,1}(M,\R)$,
called {\bf the nef cone}. A {\bf face} of the 
K\"ahler cone is the intersection
of the boundary of $\overline\Kah$ and a hyperplane $V\subset H^{1,1}(M,\R)$
which has non-empty interior.

\hfill

\conjecture (Morrison-Kawamata cone conjecture) \\
Let $M$ be a Calabi-Yau manifold. Then 
the group $\Aut(M)$ of biholomorphic automorphisms of $M$
acts on the set of faces of $\Kah$  with finite number of orbits.

\hfill

We shall be interested in the case when the manifold $M$ is simple hyperk\"ahler
(that is, IHS).
In \cite{_AV:MBM_}, we have shown that Morrison-Kawamata cone conjecture 
holds
whenever the Beauville-Bogomolov square of primitive MBM classes is 
bounded.
This is known to be the case for deformations of punctual Hilbert schemes
of K3 surfaces and for deformations of generalized Kummer varieties.
A different proof in the similar spirit for these types of 
hyperk\"ahler
manifolds has been given by Markman and
Yoshioka in \cite{MY}, under an extra assumption that the manifolds in 
question are projective.

Let us also briefly mention that this conjecture has a birational version, 
proved for projective hyperk\"ahler manifolds by E. Markman in 
\cite{_Markman:survey_} and generalized in \cite{_AV:MBM_} to the non-projective case. In this birational
version, the nef cone is replaced by the birational nef cone (that is, the
closure of the union of pullbacks of K\"ahler cones on birational models of $M$)
and the group $\Aut(M)$ is replaced by the group of birational automorphisms $\Bir(M)$. 

The key point of the proof of \cite{_AV:MBM_} is the observation that 
the orthogonal group $O(H^{1,1}_{\Z}(M), q)$ of the lattice 
$H^{1,1}_{\Z}(M)=H^{1,1}(M)\cap H^2(M,\Z)$, and therefore the 
{\bf Hodge monodromy group} $\Gamma_\Hdg$
(see \ref{monodr}) 
which is a subgroup of finite index in $O(H^{1,1}_{\Z}(M), q)$ ,
acts with finitely many orbits on the set of classes of fixed square $r\neq 0$.
When the primitive MBM classes have bounded square, we conclude that the 
monodromy
acts with finitely many orbits on the set of MBM classes. As those are 
precisely the classes whose orthogonal hyperplanes support the faces
of the K\"ahler cone, it is not difficult to deduce that there are only
finitely many, up to the action of the monodromy group, faces of the K\"ahler cone,
and also finitely many oriented faces of the K\"ahler cone (an oriented face
is a face together with the choice of normal direction). An element of the
monodromy which sends a face $F$ to a face $F'$, with both orientations pointing 
towards the 
interior of the K\"ahler cone, must preserve the K\"ahler cone.
On the other hand, Markman proved (\cite{_Markman:survey_}, Theorem 1.3)
that an element of the Hodge monodromy which
preserves the K\"ahler cone must be induced by an automorphism, so that the cone
conjecture follows.

\hfill

%\definition
%Let $M,M'$ be compact complex manifolds.
%Define a {\bf  pseudo-isomorphism} $M\dasharrow M'$ as a
%birational map which is an isomorphism outside of codimension
%$\geq 2$ subsets of $M, M'$.

%\hfill

%\remark For any pseudo-isomorphic manifolds $M,M'$,
%the second cohomologies $H^2(M)$ and $H^2(M')$ are naturally
%identified.

%\hfill

%As we have already remarked, any birational map of hyperk\"ahler
%varieties is a pseudo-isomorphism; more generally, this is true
%for all varieties with nef canonical class.

%\hfill

%\definition The

%Namely, define the
%{\bf movable cone}, also known as {\bf  birational nef cone},
%to be the closure of the union of pullbacks of $\Kah(M')$ for all
%hyperk\"ahler $M'$ birational to $M$. The Morrison-Kawamata birational cone 
%conjecture states that 
%The union of $\Kah(M')$ for all
%$M'$ pseudo-isomorphic to $M$ is called {\bf birational K\"ahler cone}.

%\hfill

%\conjecture (Morrison-Kawamata birational cone conjecture)\\
%Let $M$ be a Calabi-Yau manifold. Then 
%the group $\Bir(M)$ of birational automorphisms of $M$
%acts on the set of faces of its movable cone with finite 
%number of orbits.

%%%%%%%%%%%%%%%%%%%%%%%%%%%%%%%%%%%%%%%%%%%%%%%%%%%%%%%%%%%%
\subsection{Main results}
%%%%%%%%%%%%%%%%%%%%%%%%%%%%%%%%%%%%%%%%%%%%%%%%%%%%%%%%%%%%

The main point of the present paper is that the finiteness of the set
of primitive MBM classes of type $(1,1)$, up to the monodromy action, 
can be 
obtained without the 
boundedness assumption on their Beauville-Bogomolov square.

Our main technical result is the following 

\hfill

\theorem\label{hyperplanes} Let $L$ be a lattice of signature $(1,n)$ where $n\geq 3$,
$V=L\otimes{\mathbb R}$. Let $\Gamma$ 
be an 
arithmetic subgroup in 
$SO(1, n)$. Let $Y:= \bigcup S_i$ be a $\Gamma$-invariant union of 
rational hyperplanes $S_i$ orthogonal to negative vectors $z_i\in L$ in $V$. Then either $\Gamma$ acts on $\{S_i\}$
with finitely many orbits, or $Y$ is dense in the positive cone in $V$.

\hfill

{\bf Proof:} See \ref{_density_of_measures_SO(1,n)_Theorem_}.
\endproof

\hfill

\remark The assumption $n\geq 3$ is important for our argument which is based on Ratner theory.
We shall see that Ratner theory applies to our problem as soon as the connected component of the 
unity of $SO(1, n-1)$ is generated by unipotents, that is, for $n\geq 3$.

Taking $H^{1,1}_{\mathbb Z}(M)= H^{1,1}(M)\cap H^2(M,\Z)$ for $L$
and the Hodge monodromy group for $\Gamma$, we easily deduce: 

\hfill

\theorem\label{monodromy_acts} Assume that $M$ is projective, of Picard rank at least 4.
The monodromy group acts with finitely many orbits on the set of MBM classes 
which are
of type $(1,1)$.

\hfill

{\bf Proof:} See \ref{finitely_many_orbits}.
\endproof

\hfill

Note that, by a result of Huybrechts, the projectivity assumption for $M$ is equivalent
to the signature $(1,n)$ assumption for its Picard lattice $L$. 

\hfill

The boundedness results as an obvious corollary.

\hfill

\corollary\label{mbm_bounded} On a projective $M$ with Picard number at least 4,
primitive MBM classes of type $(1,1)$ have bounded 
Beauville-Bogomolov square.

\hfill

{\bf Proof:} Indeed, the monodromy acts by isometries.
\endproof

\hfill

Using the deformation invariance of MBM property, we can actually 
drop
the assumption that $M$ is projective and has Picard rank at least four. Indeed, if $M$ is a 
simple hyperk\"ahler manifold with $b_2(M)\geq 6$, we can always 
deform it to a projective 
manifold $M'$
on which all classes from $H^{1,1}_{\Z}(M)$ stay of type $(1,1)$ 
(see \ref{defoproj}). Since the square of a primitive MBM class is bounded 
on $M'$, the same is true for $M$.

The Morrison-Kawamata cone conjecture is then deduced as we have 
sketched it above, exactly
in the same way as in \cite{_AV:MBM_}.

\hfill

\theorem\label{mor-kaw} Let $M$ be a simple hyperk\"ahler manifold with $b_2(M)\geq 6$. 
The group of automorphisms $\Aut(M)$ acts with finitely many orbits
on the set of faces of the K\"ahler cone $\Kah(M)$.

\hfill

{\bf Proof:} See \ref{cone}.
\endproof

\hfill

\remark The theorem holds trivially for $M$ with $b_2(M)<5$, so that our result is valid 
as soon as $b_2(M)\neq 5$. This remaining case can probably be handled using methods of
hyperbolic geometry. The general belief is, though, that simple hyperk\"ahler manifolds 
with $b_2=5$ do not exist.

\hfill

Finally, as observed by Markman and Yoshioka, the boundedness of squares of primitive
MBM classes implies the following theorem (we thank Y. Kawamata for indicating
us the statement).

\hfill

\theorem\label{finitely-many-models} Let $M$ be a simple hyperk\"ahler manifold with $b_2(M)\geq 6$.
Then there are only finitely many simple hyperk\"ahler manifolds birational
to $M$.

\hfill

{\bf Proof:} This is just \cite{MY}, Corollary 1.5. Indeed, the classes $e$
menitioned in Conjecture 1.1 from \cite{MY} (that is, the classes generating
the extremal rays of the Mori cone on the simple hyperk\"ahler birational
models of $M$) are MBM classes in the sense of our \ref{mbm}.
\endproof

\hfill

The crucial tool for the proof of \ref{hyperplanes} is
Ratner theory. We recall this and some other relevant information from ergodic
theory in section \ref{ergodic}, after some preliminaries
on hyperk\"ahler manifolds in section \ref{prelim}. In section \ref{mozes-shah} we deduce \ref{hyperplanes}
from Mozes-Shah and Dani-Margulis theorems. Finally, in the last section we apply this to 
hyperk\"ahler manifolds and prove \ref{mor-kaw}.

\section{Preliminaries}\label{prelim}

%%%%%%%%%%%%%%%%%%%%%%%%%%%%%%%%%%%%%%%%%%%%%%%%%%%%%%%%%%%%
\subsection{Hyperk\"ahler manifolds, monodromy and MBM classes}
%%%%%%%%%%%%%%%%%%%%%%%%%%%%%%%%%%%%%%%%%%%%%%%%%%%%%%%%%%%% 

%%%%%%%%%%%%%%%%%%%%%%%%%%%%%%%%%%%%%%%%%%%%%%%%
\definition
A {\bf hyperk\"ahler manifold}
is a compact, K\"ahler, holomorphically symplectic manifold.

\hfill

%%%%%%%%%%%%%%%%%%%%%%%%%%%%%%%%%%%%%%%%%%%%%%%%
\definition
A hyperk\"ahler manifold $M$ is called
{\bf simple}, or IHS, if $\pi_1(M)=0$, $H^{2,0}(M)=\C$.

\hfill

This definition is motivated by Bogomolov's decomposition theorem:

\hfill

%%%%%%%%%%%%%%%%%%%%%%%%%%%%%%%%%%%%%%%%%%%%%%%%
\theorem \label{_Bogo_deco_Theorem_}
(\cite{_Bogomolov:decompo_})
Any hyperk\"ahler manifold admits a finite covering
which is a product of a torus and several 
simple hyperk\"ahler manifolds.
\endproof

\hfill

\remark
The Bogomolov decomposition theorem can be obtained
by applying the de Rham holonomy decomposition theorem
and Berger's classification of manifolds with special holonomy
to the Ricci-flat hyperk\"ahler metric on a compact
holomorphically symplectic K\"ahler manifold.
Then, a hyperk\"ahler manifold is simple
if and only if its hyperk\"ahler metric has maximal holonomy
group $\Hol(M)$ allowed by the hyperk\"ahler structure, that
is $\Hol(M)=Sp(n)$, where $n=\frac 1 2 \dim_\C M$.

\hfill

%%%%%%%%%%%%%%%%%%%%%%%%%%%%%%%%%%%%%%%%%%%%%%%
\remark
Further on, we shall assume that
all hyperk\"ahler manifolds we consider are simple.

\hfill

%\remark A hyperk\"ahler manifold naturally possesses a whole 2-sphere of 
%complex structures (see Section 2). We shall use the notation $(M,I)$ or
%$M_I$ to stress that a particular complex structure is chosen, and $M$
%to discuss the topological properties (or simply when there is no risque
%of confusion).

%\hfill

The Bogomolov-Beauville-Fujiki form was
defined in \cite{_Bogomolov:defo_} and 
\cite{_Beauville_},
but it is easiest to describe it using the
Fujiki theorem, proved in \cite{_Fujiki:HK_}.

\hfill

%%%%%%%%%%%%%%%%%%%%%%%%%%%%%%%%%%%%%%%%%%%%%%%%
\theorem\label{_Fujiki_Theorem_}
(Fujiki)
Let $M$ be a simple hyperk\"ahler manifold,
$\eta\in H^2(M)$, and $n=\frac 1 2 \dim M$. 
Then $\int_M \eta^{2n}=c q(\eta,\eta)^n$,
where $q$ is a primitive integral quadratic form on $H^2(M,\Z)$,
and $c>0$ a constant (depending on $M$). \endproof

\hfill

%%%%%%%%%%%%%%%%%%%%%%%%%%%%%%%%%%%%
\remark 
Fujiki formula (\ref{_Fujiki_Theorem_}) 
determines the form $q$ uniquely up to a sign.
For odd $n$, the sign is unambiguously determined as well.
For even $n$, one needs the following explicit
formula, which is due to Bogomolov and Beauville.
\begin{equation}\label{_BBF_expli_Equation_}
\begin{aligned}  \lambda q(\eta,\eta) &=
   \int_X \eta\wedge\eta  \wedge \Omega^{n-1}
   \wedge \bar \Omega^{n-1} -\\
 &-\frac {n-1}{n}\left(\int_X \eta \wedge \Omega^{n-1}\wedge \bar
   \Omega^{n}\right) \left(\int_X \eta \wedge \Omega^{n}\wedge \bar \Omega^{n-1}\right)
\end{aligned}
\end{equation}
where $\Omega$ is the holomorphic symplectic form, and 
$\lambda>0$.

\hfill

%%%%%%%%%%%%%%%%%%%%%%%%%%%%%%%%%%%%%%%%%%%%%%%%%%%
\definition\label{_posi_cone_Definition_}
A cohomology class $\eta \in H^{1,1}_\R(M)$ is called
{\bf negative} if $q(\eta,\eta)<0$, and {\bf positive}
if $q(\eta,\eta)>0$. Since the signature of $q$ on $H^{1,1}(M)$
is $(1, b_2-3)$, the set of positive vectors is disconnected.
{\bf The positive cone} $\Pos(M)$ is the connected component of the set
$\{\eta\in H^{1,1}_\R(M)\ \ |\ \ q(\eta,\eta)>0\}$ which contains
the classes of the K\"ahler forms. Using the
Cauchy-Schwarz inequality, it is easy to check that the positive
cone is convex.

\hfill

%%%%%%%%%%%%%%%%%%%%%%%%%%%%%%%%%%%%%%%%%%%%%%%%%%
%\remark
%Let $(M,I)$ be a hyperk\"ahler manifold,
%and $\phi:\; (M,I)\dashrightarrow (M,I')$ a bimeromorphic
%(also called ``birational'') map to another hyperk\"ahler
%manifold (note that by a result of Huybrechts 
%\cite{_Huybrechts:basic_}, birational hyperk\"ahler
%manifolds are deformation equivalent; that is why there is the same
%letter $M$ used for the source and the target of $\phi$).
% Since the canonical bundle of $(M,I)$ and $(M,I')$
%is trivial, $\phi$ is an isomorphism in codimension 1 (see for example
%\cite{_Huybrechts:basic_}, Lemma 2.6). 
%This allows one to identify $H^2(M,I)$ and $H^2(M,I')$.
%Clearly, this identification is compatible with the
%Hodge structure. Further on, we call $(M,I')$
%``a birational model'' for $(M,I)$, and identify
%$H^2(M)$ for all birational models.

\hfill

%%%%%%%%%%%%%%%%%%%%%%%%%%%%%%%%%%%%%%%%%%%%%%%%%%
%\remark
%A cohomology class $\eta\in H^2(M, \R)$ is called 
%{\bf positive} if $q(\eta,\eta)>0$, and 
%{\bf negative}, if $q(\eta,\eta)<0$.

%%%%%%%%%%%%%%%%%%%%%%%%%%%%%%%%%%%%%%%%%%%%%%%%%%%%%%%%
\definition
Let $M$ be a hyperk\"ahler manifold. The 
{\bf monodromy group} of $M$ is a subgroup of $GL(H^2(M,\Z))$
generated by the monodromy transforms for all Gauss-Manin local systems.

\hfill

It is often enlightening to consider this group in terms of the
mapping class group action. In the following paragraphs, we recall
this description.  
%(\ref{_monodro_defi_Remark_}). 

\hfill

\definition
Let $M$ be a compact complex manifold, and 
$\Diff_0(M)$ a connected component of its diffeomorphism group
({\bf the group of isotopies}). Denote by $\Comp$
the space of complex structures of K\"ahler type on $M$, equipped with
its structure of a Fr\'echet manifold (remark here that the set of complex structures of K\"ahler type 
is open in the space of all complex structures by Kodaira-Spencer stability theorem), and let
$\Teich:=\Comp/\Diff_0(M)$. We call 
it {\bf the Teichm\"uller space.} 

\hfill

For hyperk\"ahler manifolds, this is a finite-dimensional complex
non-Hausdorff manifold (\cite{_Catanese:moduli_}, \cite{_V:Torelli_}).

\hfill

\definition The {\bf mapping class group} is 
$\Diff(M)/\Diff_0(M)$. It
naturally acts on $\Teich$. The quotient of $\Teich$ by this action 
may be viewed as the ``moduli space'' for $M$. However, this space is too 
non-Hausdorff to be useful: any two open subsets of $\Teich/\Diff$
intersect (\cite{_Verbitsky:ergodic_}, \cite{_Verbitsky:ICM_}). 

\hfill

By a result of Huybrechts (see \cite{_Huybrechts:finiteness_}), in the hyperk\"ahler case $\Teich$ has only finitely many
connected components. Therefore, the subgroup of the mapping class
group which fixes the connected component of our chosen complex structure
is of finite index in the mapping class group. 

\hfill

\definition\label{monodr} The {\bf monodromy group}
$\Gamma$ is the image of this subgroup in $\Aut H^2(M,\Z)$. 
The {\bf Hodge monodromy group}
is the subgroup $\Gamma_\Hdg\subset \Gamma$ preserving the Hodge decomposition.

\hfill

The following theorem is crucial for the Morrison-Kawamata cone
conjecture.

\hfill

\theorem\label{arithmetic} (\cite{_V:Torelli_}, Theorem 3.5) 
The monodromy group is a finite index subgroup in $O(H^2(M, \Z), q)$
(and the Hodge monodromy is therefore an arithmetic subgroup of the orthogonal group of the 
Picard lattice).

\hfill

Next, we recall from \cite{_AV:MBM_} the definition of MBM classes.
Remark that any birational map between hyperk\"ahler manifolds 
$\phi: M\dasharrow M'$ is an isomorphism in codimension one
(in general this easily follows from the nefness of the canonical 
class) and therefore induces an isomorphism
on the second cohomology. We say that $M$ and $M'$ are {\bf birational
models} of each other.

\hfill

\definition\label{mbm}
A non-zero negative rational homology class
$z\in H^{1,1}(M)$ is called {\bf monodromy birationally minimal} (MBM)
if for some isometry $\gamma\in O(H^2(M,\Z))$ 
belonging to the monodromy group,
 $\gamma(z)^{\bot}\subset H^{1,1}(M)$ contains a face 
of the pull-back of the K\"ahler cone of one of birational
models $M'$ of $M$.

\hfill

%\remark
%This is equivalent to $z$ being proportional to a minimal
%rational curve in $M'$ (\cite[Theorem 5.10]{_AV:MBM_}). 
%This observation explains the
%term: MBM is an acronym for Minimal Birational Monodromy.

%\hfill

\remark Here the orthogonal is taken with respect to the 
Beauville-Bogomolov form. A face of $\Kah(M)$ is, by definition, of maximal dimension 
$h^{1,1}(M)-1$. So the definition of $z$ being MBM means that 
$\gamma(z)^{\bot}\cap\partial\Kah(M')$ contains an open subset of 
$\gamma(z)^{\bot}$. The MBM classes, or more precisely the rays they generate, 
are natural analogues of ``extremal rays'' from
projective geometry, up to monodromy and birational equivalence; hence the name.

\hfill

%%%%%%%%%%%%%%%%%%%%%%%%%%%%%%%%%%%%%%%%%%%%%
%\definition
%Let $(M,I)$ be a hyperk\"ahler manifold.
%A negative rational class $z\in H^{1,1}_\Q(M,I)$ is called
%{\bf divisorial} if $z=\lambda [D]$ for some effective divisor $D$ and% $\lambda\in \Q$.

%\hfill

The following theorem has been proved in \cite{_AV:MBM_}.

\hfill

\theorem 
Let $M$ be a hyperk\"ahler manifold, 
$z\in H^{1,1}(M)$ an integral cohomology class, $q(z,z)<0$,
and $M'$ a deformation of $M$ such that $z$ remains of type (1,1)
on $M'$. Assume that $z$ is 
monodromy birationally minimal on $M$. Then
$z$ is monodromy birationally minimal on $M'$.
%The property of $z\in H^{1,1}(M,I)$ being divisorial is likewise 
%deformation-invariant,
%provided that one restricts oneself to the complex structures with
%Picard number one (i.e. such that the Picard group is generated by $z$
%over $\Q$).
\endproof

\hfill

The MBM classes can be used to determine the K\"ahler cone of
$M$ explicitly.

\hfill

%%%%%%%%%%%%%%%%%%%%%%%%%%%%%%%%%%%%%%%%%%%%%%%%%%%%%%%
\theorem (\cite{_AV:MBM_})
Let $M$ be a hyperk\"ahler manifold, 
and $S\subset H^{1,1}(M)$ the set of all MBM classes of type $(1,1)$.
Consider the corresponding set of hyperplanes
$S^\bot:=\{W=z^\bot\ \ |\ \ z\in S\}$ in $H^{1,1}(M)$.
Then the K\"ahler cone of $M$ is 
a connected component of $\Pos(M)\backslash \cup S^\bot$,
where $\Pos(M)$ is the positive cone of $M$.
Moreover, for any connected component $K$ of 
$\Pos(M)\backslash \cup S^\bot$,
there exists $\gamma\in O(H^2(M,\Z))$ in the monodromy group of $M$ and
a birational model $M'$ such that $\gamma(K)$ is the K\"ahler cone 
of $M'$.
\endproof

\hfill

\remark The main point of this theorem is that for a negative integral
class $z\in H^{1,1}(M)$, the orthogonal hyperplane either passes
through the interiour of some K\"ahler-Weyl chamber and then it contains no face of a K\"ahler-Weyl chamber 
(that is, $z$ is not
MBM), or its 
intersection with the positive cone is a union of faces of such 
chambers (when $z$ is MBM). This is illustrated by a picture
taken from \cite{_AV:MBM_}:

\centerline{\begin{tabular}{cc}
\includegraphics[width=0.30\textwidth]{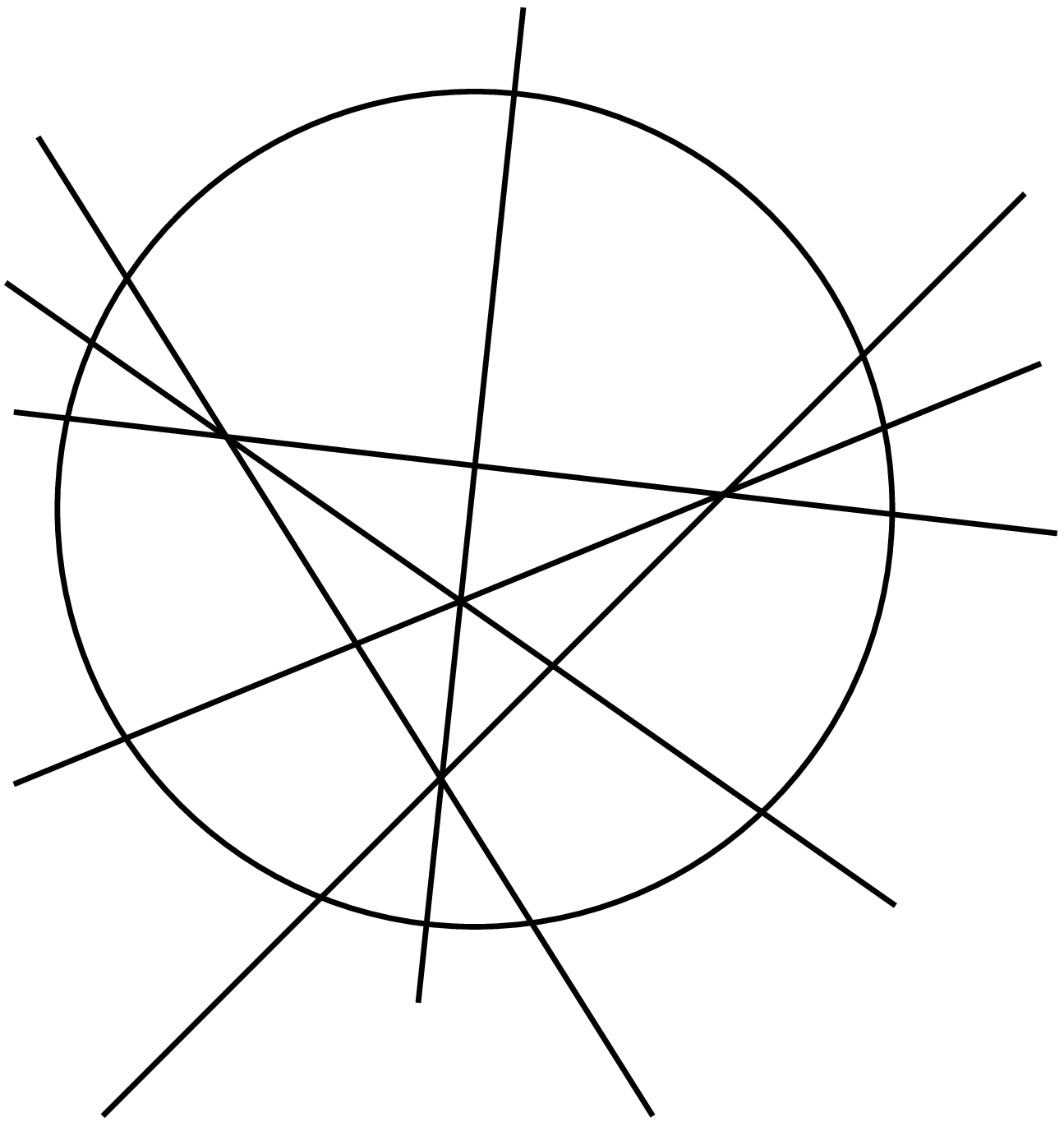}\  
&\  \includegraphics[width=0.30\textwidth]{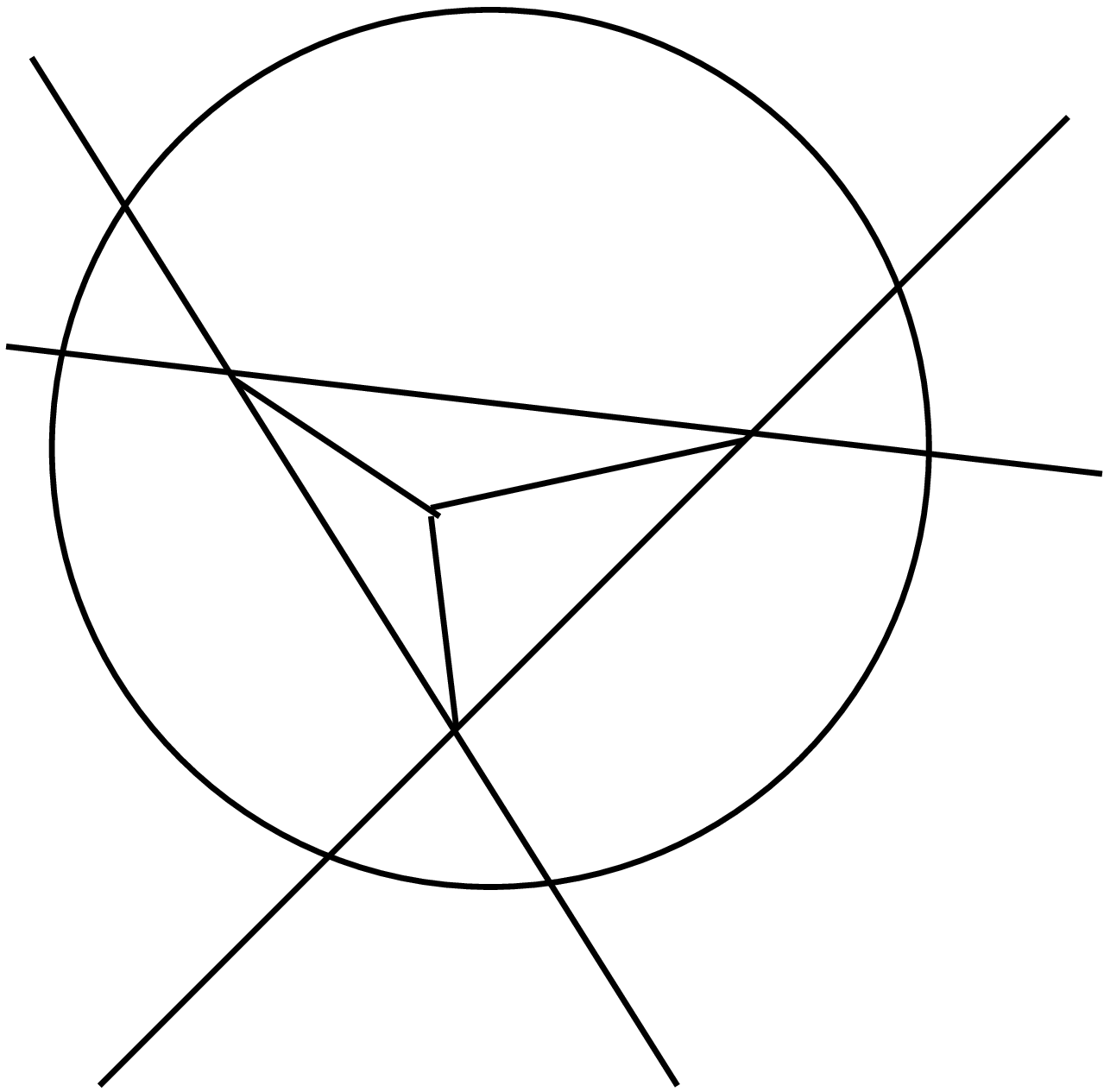} \\
{\bf \scriptsize Allowed partition} & {\bf \scriptsize Prohibited partition} 
\end{tabular}}

%%%%%%%%%%%%%%%%%%%%%%%%%%%%%%%%%%%%%%%%%%%%%%%%%%%%%%%%%%%%%%%%%%%%%%%%

\subsection{Global Torelli theorem and deformations}

%%%%%%%%%%%%%%%%%%%%%%%%%%%%%%%%%%%%%%%%%%%%%%%%%%%%%%%%%%%%%%%%%%%%%%%%

In this subsection, we recall a number of 
results  about deformations of hyperk\"ahler manifolds
used further on in this paper. For more
details and references,  see \cite{_V:Torelli_}.

\hfill

%%%%%%%%%%%%%%%%%%%%%%%%%%%%%%%%%%%%%%%%%%%%%%%%%%%%%%%%%%%%
%\definition
%Let $M$ be a compact complex manifold, and 
%$\Diff_0(M)$ a connected component of its diffeomorphism group
%({\bf the group of isotopies}). Denote by $\Comp$
%the space of complex structures on $M$, equipped with
%its structure of a Fr\'echet manifold, and let
%$\Teich:=\Comp/\Diff_0(M)$. We call 
%it {\bf the Teichm\"uller space.} 

%\hfill

%\remark
%In many important cases, such as
%for manifolds with trivial canonical class (\cite{_Catanese:moduli_}), 
%$\Teich$ is a finite-dimensional
%complex space; usually it is non-Hausdorff.

%\hfill

%\definition\label{univfam} The universal family
%of complex manifolds over $\Teich$ is defined as $\Univ/\Diff_0(M)$,
%where $\Univ$ is the natural universal family over $\Comp$ with its 
%Fr\'echet manifold structure. Locally, it is isomorphic to the universal
%family over the Kuranishi 
%space.\footnote{We are grateful to Claire Voisin for this observation.}

%\hfill

%\definition The {\bf mapping class group} is $\Gamma=\Diff(M)/\Diff_0(M)$. It
%naturally acts on $\Teich$ (the quotient of $\Teich$ by $\Gamma$ may be
%viewed as the ``moduli space'' for $M$, but in general it has very bad 
%properties; see below).

%\hfill

Let $M$ be a hyperk\"ahler manifold (as usual, we assume
$M$ to be simple). Any deformation $M'$ of $M$
is also a simple hyperk\"ahler manifold, 
because the Hodge numbers are constant in families and 
thus $H^{2,0}(M')$ is one-dimensional. Let us view $M'$ as a couple 
$(M, J)$, where $J$ is a new complex structure on $M$, 
that is, a point of the Teichm\"uller space 
$\Teich$.

\hfill

%%%%%%%%%%%%%%%%%%%%%%%%%%%%%%%%%%%%%%%%%%%%%%%%
\definition
 Let 
\[ \Per:\; \Teich \arrow {\Bbb P}H^2(M, \C)
\]
map $J$ to the line $H^{2,0}(M,J)\in {\Bbb P}H^2(M, \C)$.
The map $\Per$ is 
called {\bf the period map}.

\hfill

%%%%%%%%%%%%%%%%%%%%%%%%%%%%%%%%%%%%%%%%%%%%%%%%
\remark
The period map $\Per$ maps $\Teich$ into an open subset of a 
quadric, defined by
\[
\Perspace:= \{l\in {\Bbb P}H^2(M, \C)\ \ | \ \  q(l,l)=0, q(l, \bar l) >0\}.
\]
It is called {\bf the period domain} of $M$.
Indeed, any holomorphic symplectic form $l$
satisfies the relations $q(l,l)=0, q(l, \bar l) >0$,
as follows from \eqref{_BBF_expli_Equation_}.

\hfill

%%%%%%%%%%%%%%%%%%%%%%%%%%%%%%%%%%%%%%%%%%%%%%%%
%\proposition\label{_period_Grassmann_Proposition_}
%The period domain $\Perspace$
%is identified with the quotient
%$SO(b_2-3,3)/SO(2) \times SO(b_2-3,1)$, which
%is the Grassmannian of positive oriented 2-planes in $H^2(M,\R)$.

%\hfill

%{\bf Proof:} This statement is well known, but we shall
%sketch its proof for the reader's convenience.

%{\bf Step 1:} Given $l\in {\Bbb P}H^2(M, \C)$, the space
%generated by $\Im l, \Re l$ is 2-dimensional, because 
%$q(l,l)=0, q(l, \bar l)>0$ implies that $l \cap H^2(M,\R)=0$.

%{\bf Step 2:}  This 2-dimensional plane is 
%positive, because 
% $q(\Re l, \Re l) = q(l+ \bar l, l+ \bar l) = 2 q(l, \bar l)>0$.

%{\bf Step 3:} Conversely, for any 2-dimensional positive
%plane  $V\in H^2(M,\R)$, 
%the quadric $\{l\in V \otimes_\R \C\ \ | \ \ q(l,l)=0\}$
%consists of two lines; a choice of a line is determined by the
%orientation.
%\endproof

\hfill

\definition
Let $M$ be a topological space. We say that $x, y \in M$
are {\bf non-separable} (denoted by $x\sim y$)
if for any open sets $V\ni x, U\ni y$, $U \cap V\neq \emptyset$.

\hfill

By a result of Huybrechts \cite{_Huybrechts:basic_}, non-separable points of $\Teich$ correspond to birational 
hyperk\"ahler manifolds.

%%%%%%%%%%%%%%%%%%%%%%%%%%%%%%%%%%%%%%%%%%%%%%%%%%%%%%%%%%%%
%\theorem\label{nonsep-birat} 
%(Huybrechts; \cite{_Huybrechts:basic_}).
%If two points $I,I'\in \Teich$ are non-separable, then  
%there exists a bimeromorphism $(M,I)\dasharrow (M,I')$.
%\endproof

%\hfill

%\remark\label{birat-nonsep}
%The converse is not true since many points in the Teichm\"uller 
%space correspond to isomorphic manifolds, and they are not always pairwise
%non-separable (consider orbits of the
%mapping class group action). But it is true that if 
%two hyperk\"ahler manifolds are bimeromorphic, then one can find 
%non-separable points in the Teichm\"uller space representing them.

\hfill

\definition
The space $\Teich_b:= \Teich\!/\!\!\sim$ is called {\bf the
birational Teichm\"uller space} of $M$.

\hfill

\remark This terminology is slightly misleading since there are 
non-separable points of the Teichm\"uller space which correspond
to biregular, not just birational, complex structures. Even for
K3 surfaces, the Teichm\"uller space is non-Hausdorff. 

\hfill

%%%%%%%%%%%%%%%%%%%%%%%%%%%%%%%%%%%%%%%%%%%%%%%%%%%%%%%%%%%%
\theorem \label{_glo_Torelli_Theorem_}
(Global Torelli theorem; \cite{_V:Torelli_})
The period map 
$\Teich_b\stackrel \Per \arrow \Perspace$ is an isomorphism
on each connected component of $\Teich_b$.
\endproof

\hfill

 By a result of Huybrechts (\cite{_Huybrechts:finiteness_}), $\Teich$ has only finitely
many connected components. We shall fix the component  ${\Teich}^0$ 
containing the parameter point for our initial complex structure, and denote by 
$\tilde{\Gamma}$
the subgroup of finite index in the mapping class group fixing this component.

It is natural to view the quotient of $\Teich$ by the mapping class group as a
moduli space for $M$ and the quotient of $\Teich_b$  by the mapping class group as a ``birational moduli
space'': indeed its points are in bijective correspondence with the complex structures of hyperk\"ahler type on $M$
up to a bimeromorphic equivalence.

\hfill

%%%%%%%%%%%%%%%%%%%%%%%%%%%%%%%%%%%%%%%%%%%%%%%%%%%%%%%%%%%%
%\definition
%Let $M$ be a hyperkaehler manifold,
%$\Teich_b$ its birational Teichm\"uller space,
%and $\Gamma$ the {\bf mapping class group} $\Diff (M)/\Diff_0(M)$.
%The quotient $\Teich_b/\Gamma$ is called
%{\bf the birational moduli space} of $M$.
%Its points are in bijective correspondence with the
%complex structures of hyperk\"ahler type on $M$
%up to a bimeromorphic equivalence.

%\hfill

%%%%%%%%%%%%%%%%%%%%%%%%%%%%%%%%%%%%%%%%%%%%%%%%
\remark
The word ``space'' in this context is misleading.
In fact, the quotient topology on $\Teich^0_b/\tilde{\Gamma}$ is extremely
non-Hausdorff, e.g. every two open sets would intersect
(\cite{_Verbitsky:ergodic_}).

\hfill

The Global Torelli theorem can be stated
as a result about the birational moduli space.

\hfill

%%%%%%%%%%%%%%%%%%%%%%%%%%%%%%%%%%%%%%%%%%%%%%%%
\theorem\label{_moduli_monodro_Theorem_}
(\cite[Theorem 7.2, Remark 7.4, Theorem 3.5]{_V:Torelli_})
Let $(M,I)$ be a hyperk\"ahler manifold, and $W$ 
a connected component of its birational
moduli space. Then $W$ is isomorphic to ${\Perspace}/\Gamma$,
where $\Gamma$ is 
an arithmetic subgroup in $O(H^2(M, \R), q)$, called {\bf the
monodromy group} of $(M,I)$. In fact $\Gamma$ is the image of $\tilde{\Gamma}$
in $O(H^2(M, \R), q)$.  
\endproof

\hfill

%%%%%%%%%%%%%%%%%%%%%%%%%%%%%%%%%%%%%%%%%%%%%%%%
\remark\label{_monodro_defi_Remark_} As we have already mentioned,
the monodromy group of $(M,I)$ can be also described
as a subgroup of the group $O(H^2(M, \Z), q)$
generated by monodromy transform maps for 
Gauss-Manin local systems obtained from all
deformations of $(M,I)$ over a complex base
(\cite[Definition 7.1]{_V:Torelli_}). This is 
how this group was originally defined by Markman
(\cite{_Markman:constra_}, \cite{_Markman:survey_}).
%The fact that it is of finite index in $O(H^2(M, \Z), q)$
%is crucial for the Morrison-Kawamata conjecture, see next section. 

\hfill

%%%%%%%%%%%%%%%%%%%%%%%%%%%%%%%%%%%%%%%%%%%%%%%%
%\remark 
%A caution: usually ``the global Torelli theorem''
%is understood as a theorem about Hodge structures.
%For K3 surfaces, the Hodge structure on $H^2(M,\Z)$
%determines the complex structure. 
%For $\dim_\C M >2$, it is false.

\definition Let $z\in H^2(M,\Z)$ be an integral cohomology class. The space
$\Teich_z$ is the part of $\Teich$ where the class $z$ is of type $(1,1)$.

\hfill

The following proposition is well-known.

\hfill

\proposition\label{teichz}
$\Teich_z$ is the inverse image under the period map of the subset ${\Perspace}_z\subset \Perspace$
which consists of $l$ with $q(l,z)=0$.

\hfill 

{\bf Proof:} This is clear since $H^{1,1}(M)$ is the orthogonal, under $q$, to $H^{2,0}(M)\oplus H^{0,2}(M)$.
\endproof

\hfill

By a theorem of Huybrechts, a holomorphic symplectic manifold $M$ is projective if and only if it has an
integral $(1,1)$-class with strictly positive Beauville-Bogomolov square. In this case, the Picard lattice 
$H^{1,1}_{\Z}(M)=
H^2(M, \Z)\cap H^{1,1}(M)$, equipped with the Beauville-Bogomolov form $q$, is a lattice of signature
$(+, -,-,\dots, -)$. If $M$ is not projective, the Picard lattice can be either negative definite, or
degenerate negative semidefinite with one-dimensional kernel. In both cases, its rank cannot be maximal
(i.e. equal to the dimension of $H^{1,1}(M)$), since the signature of $q$ on $H^{1,1}(M)$ is $(+, -,-,\dots, -)$.
Together with this observation, \ref{teichz} easily implies the following

\hfill

\proposition\label{defoproj}
Let $M$ be an irreducible holomorphic symplectic manifold. There exists a deformation $M'$ of $M$ which is
projective and such that all integral $(1,1)$-classes on $M$ remain of type $(1,1)$ on $M'$.
Moreover one can take $M'$ of maximal Picard rank $h^{1,1}(M)$.

\hfill 

{\bf Proof:} By \ref{teichz}, the locus $C$ where all integral $(1,1)$-classes on $M$ remain of type 
$(1,1)$ is the preimage of the intersection of $N$ complex hyperplanes and ${\Perspace}$, where $N$ is strictly less
than the (complex) dimension of ${\Perspace}$. It is therefore strictly positive-dimensional. For $M'$ representing
a general point of
this locus, the Picard lattice is the same as that of $M$, but at a special point the Picard number jumps.
Namely it jumps along the intersection with each hyperplane of the form $z^{\perp}$, where $z$ is an integral
$(1,1)$-class. In particular, there are isolated points inside $C$ where the Picard rank is maximal.
By the observations above, the corresponding variety $M'$ must be projective.
\endproof 

\hfill

This proposition shall be useful in reducing the cone conjecture to the projective case with 
high Picard number (see \ref{boundedness-nonproj}).

%%%%%%%%%%%%%%%%%%%%%%%%%%%%%%%%%%%%%%%%%%%%%%%%%%%%%%%%%%%%
\section{Ergodic theory and its applications}\label{ergodic}
%%%%%%%%%%%%%%%%%%%%%%%%%%%%%%%%%%%%%%%%%%%%%%%%%%%%%%%%%%%% 

%%%%%%%%%%%%%%%%%%%%%%%%%%%%%%%%%%%%%%%%%%%%%%%%%%%%%%%%%%%%
\subsection{Ergodic theory: basic definitions and facts}
%%%%%%%%%%%%%%%%%%%%%%%%%%%%%%%%%%%%%%%%%%%%%%%%%%%%%%%%%%%% 

%%%%%%%%%%%%%%%%%%%%%%%%%%%%%%%%%%%%%%%%%%%%%%%%%%
\definition
Let $(M,\mu)$ be a space with a measure,
and $G$ a group acting on $M$ preserving $\mu$.
This action is {\bf ergodic} if all
$G$-invariant measurable subsets $M'\subset M$
satisfy $\mu(M')=0$ or $\mu(M\backslash M')=0$.

\hfill

%%%%%%%%%%%%%%%%%%%%%%%%%%%%%%%%%%%%%%%%%%%%%%%%%%%%%%%%%%%%
\claim\label{_non_dense_m_zero_Claim_}
Let $M$ be a manifold, $\mu$ the Lebesgue measure, and
$G$ a group acting on $(M,\mu)$ ergodically.  Then the 
set of points with non-dense orbits has measure 0.

\hfill

{\bf Proof:}
Consider a non-empty open subset $U\subset M$. 
Then $\mu(U)>0$, hence $M':=G\cdot U$ satisfies 
$\mu(M\backslash M')=0$. For any orbit $G\cdot x$
not intersecting $U$, $x\in M\backslash M'$.
Therefore the set of such points has measure 0.
\endproof

\hfill

%%%%%%%%%%%%%%%%%%%%%%%%%%%%%%%%%%%%%%%%%%%%
\definition 
Let $G$ be a Lie group, and $\Gamma\subset G$ a discrete
subgroup. Consider the pushforward of the Haar measure to
$G/\Gamma$. Here, by abuse of terminology, 
``taking the pushforward'' of a measure
means measuring the intersection of the 
inverse image with a fixed fundamental 
domain. We say that $\Gamma$ {\bf has finite covolume}
if the Haar measure of $G/\Gamma$ is finite.
In this case $\Gamma$ is called {\bf a lattice subgroup}.

\hfill

%%%%%%%%%%%%%%%%%%%%%%%%%%%%%%%%%%%%%%%%%
\remark\label{_Borel_HCh_Remark_}
Borel and Harish-Chandra proved that
an arithmetic subgroup of a reductive group $G$
is a lattice whenever $G$ has no non-trivial characters
over $\Q$ (see \cite{bo-hch}, Theorem 7.8 for semisimple case and Theorem 9.4 for the general case).
%e.g. \cite{_Vinberg_Gorbatsevich_Shvartsman_}). 
In particular, all arithmetic subgroups
of a semi-simple group are lattices. Therefore the monodromy and the Hodge monodromy groups
from the previous section
are lattices in the corresponding orthogonal groups, which is a very important point for us.

\hfill

In this paper, we deal with the following example of an ergodic action.

\hfill

%%%%%%%%%%%%%%%%%%%%%%%%%%%%%%%%%%%%%%%%%%%%%%%%
\theorem \label{_Moore_Theorem_}
(Calvin C. Moore, \cite[Theorem 4]{_Moore:ergodi_})
Let $\Gamma$ be a lattice subgroup 
(such as an arithmetic subgroup) in a non-compact 
simple Lie group $G$ with finite center, and $H\subset G$ a 
Lie subgroup. Then the left action of $H$
on $G/\Gamma$ is ergodic if and only if the closure of $H$ is non-compact. \endproof

\hfill

Let us also state the following classical result.

\hfill

\theorem \label{birkhoff} (Birkhoff ergodic theorem, see for example \cite{W}, 
1.6) Let $\mu$ be a probability
measure on a 
%locally compact and separated 
manifold $X$, and let $g_t$
be an ergodic flow preserving $\mu$. Then for almost all $x\in X$ and any
$f\in L^1(\mu)$, the
limit of $m_T(f)=\frac {1}{T}\int^T_0 f(g_tx)dt$ as $T\to +\infty$  exists and equals 
$\int_Xfd\mu$. In particular, for any measurable subset $K$ and almost all $x$, 
the part
of time that the orbit of $x$ spends in $K$ is equal to $\mu(K)$.
\endproof

\hfill

%%%%%%%%%%%%%%%%%%%%%%%%%%%%%%%%%%%%%%%%%%%%%%%%
\subsection{Lie groups generated by unipotents}
%%%%%%%%%%%%%%%%%%%%%%%%%%%%%%%%%%%%%%%%%%%%%%%%

Here we state some of the main results of Ratner theory.
We follow \cite{_Kleinbock_etc:Handbook_} 
and \cite{_Morris:Ratner_}.

\hfill

%%%%%%%%%%%%%%%%%%%%%%%%%%%%%%%%%%%%%%%%%%%%%%%%
\definition
Let $G$ be a Lie group, and $g\in G$ any element.
We say that $g$ is {\bf unipotent} if $g=e^h$ for a
nilpotent element $h$ in its Lie algebra.
A group $G$ is {\bf generated by unipotents}
if $G$ is multiplicatively generated by unipotent one-parameter subgroups.

\hfill

%%%%%%%%%%%%%%%%%%%%%%%%%%%%%%%%%%%%%%%%%%%%%%%%%%%%%
\theorem\label{_ratner_orbit_clo_Theorem_}
(Ratner orbit closure theorem, \cite{R1})\\
Let $H\subset G$ be a Lie subroup generated by 
unipotents, and $\Gamma\subset G$ a lattice.
Then the closure of any $H$-orbit $Hx$ in $G/\Gamma$
is an orbit of a closed, connected subgroup $S\subset G$,
such that $S\cap x\Gamma x^{-1}\subset S$ is a lattice in $S$.

\hfill

{\bf Proof:} \cite[1.1.15 (2)]{_Morris:Ratner_}. \endproof

\hfill

For arithmetic groups
Ratner orbit closure theorem can be stated in a more
precise way, as follows.

\hfill

%%%%%%%%%%%%%%%%%%%%%%%%%%%%%%%%%%%%%%%%%%%%%%%%
\theorem\label{_closure_arithm_Ratner_Theorem_}
Let $G$ be a real algebraic group defined over $\Q$
and with no non-trivial characters, $W\subset G$ a
subgroup generated by unipotents, and $\Gamma\subset G$
an arithmetic lattice. For a given $g\in G$,
let $H$ be the smallest real algebraic $\Q$-subgroup of
$G$ containing $g^{-1}Wg$. Then the closure of 
$Wg$ in $G/\Gamma$ is $Hg$. 

\hfill

{\bf Proof:} See \cite[Proposition 3.3.7]{_Kleinbock_etc:Handbook_}
or \cite[Proposition 3.2]{_Shah:uniformly_}.
\endproof

\hfill

Ratner orbit closure theorem is a consequence of her fundamental result on ergodic measures \cite{R2},
known as Ratner measure classification theorem, which we recall below.

\hfill

\definition\label{_algebraic_measure_Definition_}
Let $G$ be a Lie group, $\Gamma$ a lattice, and $G/\Gamma$
the quotient space, considered as a space with Haar measure.
Consider an orbit $S\cdot x\subset G$ of a closed subgroup $S\subset G$, 
put the Haar measure on $S\cdot x$, and assume that its image
in $G/\Gamma$ has finite Haar measure (this means that
$S\cap x\Gamma x^{-1}$ is a lattice in $S$).
A measure on $G/\Gamma$ is called {\bf algebraic}
if it is proportional to the pushforward of the Haar
measure on $S\cdot x/\Gamma$ to
$G/\Gamma$. 

\hfill

If $G$ is a non-compact simple Lie group with 
finite center and $H\subset G$ is a Lie subgroup
with non-compact closure, as in Moore's theorem 
(\ref{_Moore_Theorem_}), consider the algebraic
measure on $G/\Gamma$ which is proportional to 
the pushforward of the Haar measure of $S$, 
where $S$ is taken from the Ratner's orbit closure theorem. 
It follows from Moore's theorem
that the action of $H$ on $G/\Gamma$ is ergodic. 
Ratner's measure classification theorem states 
that all invariant ergodic measures  under 
the action of subgroups {\bf generated by unipotents} arise in this way. 

\hfill

%%%%%%%%%%%%%%%%%%%%%%%%%%%%%%%%%%%%%%%%%%%%%%%%%%%%%%%%%%%%
\theorem\label{ratner-ergodic}
(Ratner's measure classification theorem, \cite{R2})\\
Let $G$ be a connected Lie group, $\Gamma$ a lattice, and $G/\Gamma$
the quotient space, considered as a space with Haar measure.
Consider a finite measure $\mu$ on $G/\Gamma$.
Assume that $\mu$ is invariant and ergodic with respect to an action 
of a subgroup $H\subset G$ generated by unipotents.
Then $\mu$ is algebraic.

\hfill

{\bf Proof:} see \cite[1.3.7]{_Morris:Ratner_}. \endproof

\hfill

\remark\label{monoparam}
In most texts, Ratner theorems are formulated for {\bf unipotent flows}, that is, $H$
is assumed to be a one-parameter unipotent subgroup $\{u(t)|t\in\R \}$. One gets rid of this 
assumption using the following lemma.

\hfill

\lemma\label{_one-param_Lemma_}
(\cite[Lemma 2.3]{_Mozes_Shah_} or 
\cite[Corollary 3.3.5]{_Kleinbock_etc:Handbook_}) 
Let $H$ be a subgroup of $G$ generated by unipotent one-parameter subgroups. Then any finite
$H$-invariant $H$-ergodic measure is ergodic with respect to some one-parameter unipotent subgroup
of $H$.

%%%%%%%%%%%%%%%%%%%%%%%%%%%%%%%%%%%%%%%%%%%%%%%%%%%%%%%%%%%%

\section{Algebraic measures on homogeneous spaces}\label{mozes-shah}

%%%%%%%%%%%%%%%%%%%%%%%%%%%%%%%%%%%%%%%%%%%%%%%%%%%%%%%%%%%%

The main result of this section
(\ref{_density_of_measures_SO(1,n)_Theorem_}) follows from a theorem of
Mozes and Shah \cite[Theorem 1.1]{_Mozes_Shah_}. 
%A more relevant version, but for a different (less general and more complicated) situation
%is found in \cite[Theorem 2.3]{_Eskin_Mirzakhani_Mohammadi_}.

%%%%%%%%%%%%%%%%%%%%%%%%%%%%%%%%%%%%%%%%%%%%%%%%%%%%%%%%%%%%
\subsection{Limits of ergodic measures}
%%%%%%%%%%%%%%%%%%%%%%%%%%%%%%%%%%%%%%%%%%%%%%%%%%%%%%%%%%%%

%%%%%%%%%%%%%%%%%%%%%%%%%%%%%%%%%%%%%%%%%%%%%%%%%%%%%%%%%%%%
%\definition\label{_algebraic_measure_Definition_}
%Let $G$ be a Lie group, $\Gamma$ a lattice, and $G/\Gamma$
%the quotient space, considered as a space with Haar measure.
%Consider an orbit $S\cdot x\subset G$ of a subgroup $S\subset G$, 
%put the Haar measure on $S\cdot x$, and assume that its image
%in $G/\Gamma$ has finite Haar measure (this means that
%$S\cap x\Gamma x^{-1}$ is a lattice in $S$).
%A measure on $G/\Gamma$ is called {\bf algebraic},
%if it is proportional to the pushforward of the Haar
%measure on $S\cdot x/\Gamma$ to
%$G/\Gamma$.

%\hfill

%%%%%%%%%%%%%%%%%%%%%%%%%%%%%%%%%%%%%%%%%%%%%%%%%%%%%%%%%%%%
%\theorem
%(Ratner's ergodic theorem)\\
%Let $G$ be a Lie group, $\Gamma$ a lattice, and $G/\Gamma$
%the quotient space, considered as a space with Haar measure.
%Consider a finite measure $\mu$ on $G/\Gamma$.
%Assume that $\mu$ is ergodic with respect to an action 
%of a subgroup $H\subset G$ generated by unipotents.
%Then $\mu$ is algebraic. Moreover, any algebraic
%measure associated with a subgroup $S\subset G$ 
%as in \ref{_algebraic_measure_Definition_} 
%is ergodic, if $S$ is generated by unipotents.

%\hfill

%{\bf Proof:} \cite[???]???. \endproof

\definition
Recall that {\bf a Polish topological space}
is a metrizable topological space with countable base.
Let $V$ be the set of all finite Borel measures on a 
Polish topological space $M$,
and $C^0(M)$ the space of bounded continuous functions.
{\bf Weak topology} on $V$ is the weakest topology 
in which for all $f\in C^0(M)$ the
map $V \arrow \R$ given by $\mu \arrow \int_M f\mu$ is 
continuous. If one identifies $V$ with a subset 
in $C^0(M)^*$, the weak topology is identified with the weak-*
topology on $C^0(M)^*$. This is why it is also called 
{\bf the weak-* topology}. 

\hfill

\remark
It is not hard to prove that
the space of probability measures on a compact Polish space is
compact in weak topology. This explains the usefulness of this notion.

\hfill

%%%%%%%%%%%%%%%%%%%%%%%%%%%%%%%%%%%%%%%%%%%%%%%%%%%%%%%%%%%%
\theorem\label{_MS_main_Theorem_}
(Mozes-Shah theorem)\\
Let $G$ be a connected Lie group, $\Gamma$ a lattice, 
 $\{u_i(t)\}\subset G$ a sequence of unipotent
one-parameter subgroups in $G$, and $\mu_i$ a sequence
of $u_i$-invariant, $u_i$-ergodic probability measures on
$G/\Gamma$, associated with orbits $S_i\cdot x_i\subset G/\Gamma$
as in \ref{_algebraic_measure_Definition_}. 
Assume that $\lim \mu_i=\mu$
with respect to weak topology, with $\mu$ a probability measure on $X$, and let $x\in \Supp(\mu)$. Then
\begin{description}
\item[(i)] $\mu$ is an algebraic measure, associated with an orbit $S\cdot x$
as in \ref{_algebraic_measure_Definition_}. 

\item[(ii)] Let $g_i\in G$ be elements which satisfy $g_i x_i=x$, and assume that $g_i\to e$ in $G$
(so that $x_i$ converge to $x$).
Then there exists $i_0\in {\Bbb N}$ such that for all
$i>i_0$, $S\cdot x \supset g_i S_i \cdot x_i$.
\end{description}

{\bf Proof:} The statement (i) follows from 
\cite[Theorem 1.1 (3)]{_Mozes_Shah_} and Ratner
measure classification theorem, and (ii) is \cite[Theorem 1.1 (2)]{_Mozes_Shah_}.
\endproof

\hfill

\remark More precisely, in \cite[Theorem 1.1]{_Mozes_Shah_} there is an additional condition that
the trajectories $\{u_i(t)\}x_i,\ t>0$ should be uniformly distributed with respect to $\mu_i$.
But this is automatic by another theorem of Ratner (Ratner equidistribution theorem, see e.g. 
\cite{_Morris:Ratner_}, Theorem 1.3.4), and in fact already by Birkhoff ergodic
theorem (\ref{birkhoff}), which states the uniform distribution of orbits of one-parameter
subgroups for almost all starting points.

\hfill

%%%%%%%%%%%%%%%%%%%%%%%%%%%%%%%%%%%%%%%%%%%%%%%%%%%%%%%%%%%%

The following theorem is an interpretation of Dani-Margulis theorem
as stated in  \cite[Theorem 6.1]{_Dani_Margulis_} obtained by applying 
Birkhoff ergodic therorem.

\hfill

\theorem\label{_Dani_Margulis_Theorem_}
(Dani-Margulis theorem). \\
 Let $G$ be a connected Lie group, $\Gamma$ a lattice,
$X:= G/\Gamma$, $C\subset X$ a compact subset,
and $\epsilon >0$. Then there exists a
compact subset $K\subset X$ such that 
for any algebraic probability measure
$\mu$ on $X$, satisfying $\mu(C)\neq 0$
and associated with a group generated by unipotents,
one has $\mu(K)\geq 1-\epsilon$.
\endproof

\hfill

{\bf Proof:} By \ref{_one-param_Lemma_}, 
$\mu$ is invariant and ergodic with respect
to a one-parameter unipotent subgroup $u(t)$. Now apply
\cite[Theorem 6.1]{_Dani_Margulis_} to a starting point $x$ which is 
one of ``almost all points'' of $\Supp (\mu)\cap C$ in the sense of Birkhoff
theorem. \endproof

\hfill

Combining Dani-Margulis theorem and Mozes-Shah theorem, one gets the following useful corollary
(\cite[Corollary 1.1, Corollary 1.3, Corollary 1.4]{_Mozes_Shah_}).

\hfill

%%%%%%%%%%%%%%%%%%%%%%%%%%%%%%%%%%%%%%%%%%%%%%%%%%%%%%%%%%%%
\corollary\label{_compactness_MS_Theorem_}
Let $G$ be a connected Lie group, $\Gamma$ a lattice,
${\cal P}(X)$ be the space of all probability measures
on $X=G/\Gamma$, and ${\cal Q}(X)\subset {\cal P}(X)$ the space of all
algebraic probability measures associated with subgroups $H\subset G$
generated by unipotents (as in Ratner theorems). Then ${\cal Q}(X)$
is closed in ${\cal P}$ with respect to weak-star 
topology. 
%and the set ${\cal Q}(x)$ of all $\mu\in {\cal Q}$
%which satisfy $\Supp(\mu)\ni x$ is compact. 
Moreover, let $X\cap \{\infty\}$
denote the one-point compactification of $X$, so that ${\cal P}(X\cap \{\infty\})$ is compact.
If for a sequence $\mu_i\in {\cal Q}(X)$, $\mu_i\to \mu \in {\cal P}(X\cap \{\infty\})$, then 
either $\mu \in {\cal Q}(X)$,
or $\mu$ is supported at infinity. 

\endproof

%%%%%%%%%%%%%%%%%%%%%%%%%%%%%%%%%%%%%%%%%%%%%%%%%%%%%%%%%%%%
\subsection{Rational hyperplanes intersecting a compact set}
%%%%%%%%%%%%%%%%%%%%%%%%%%%%%%%%%%%%%%%%%%%%%%%%%%%%%%%%%%%%

%%%%%%%%%%%%%%%%%%%%%%%%%%%%%%%%%%%%%%%%%%%%%
\definition
Let $V_\Q$ be an $n+1$-dimensional 
rational vector space with a scalar product
of signature $(+, -,-, ..., -)$, and $V:= V_\Q\otimes_\Q \R$.
We consider the projectivization of the positive cone ${\Bbb P}^+ V$
as the hyperbolic space of dimension $n$. Given a $k+1$-dimensional
subspace $W_\Q\subset V_\Q$ such that the restriction of the scalar product to $W_\Q$ 
still has signature
$(1,k)$, we may associate the projectivized
positive cone ${\Bbb P}^+ W\subset{\Bbb P}^+ V$ with $W=W_\Q \otimes_\Q \R$. 
When $k=n-1$, we shall call ${\Bbb P}^+ W\subset{\Bbb P}^+ V$
{\bf a rational hyperplane} in ${\Bbb P}^+ V$.

\hfill

Let $\Gamma$ be a rational lattice in the group of isometries
of ${\Bbb P}^+ V$, and $\{S_i\}$ a set of rational hyperplanes.
We are interested in the images of $S_i$ in ${\Bbb P}^+ V/\Gamma$.
The following theorem can be used to show that these images all
intersect a compact subset of ${\Bbb P}^+ V/\Gamma$.

\hfill

%%%%%%%%%%%%%%%%%%%%%%%%%%%%%%%%%%%%%%%%%%%%%%%%
\theorem\label{_intersect_compact_Theorem_}
Let $\{S_i\}$ be a set of rational hyperplanes in ${\Bbb P}^+ V$,
$P_\Q\subset V_\Q$ a rational subspace of signature $(1,2)$, and
${\Bbb P}^+ P\subset {\Bbb P}^+ V$ the corresponding 2-dimensional
hyperbolic subspace. Consider an arithmetic lattice 
$\Gamma\subset SO(V, \Z)$, and let $\Gamma_P$ be the stabilizer
of $P_\Q$ in $\Gamma$. Then there exists a compact subset
$K\subset {\Bbb P}^+ P$ such that $\Gamma_P\cdot K$ intersects
all the hyperplanes $S_i$.

\hfill

{\bf Proof:} Since $\Gamma$ has finite index in a lattice $O(V,\Z)$,
$\Gamma_P$ has finite index in the lattice $O(P,\Z)$. One may view $\Gamma_P$ as a
multi-dimensional analogue of Fuchsian or Kleinian
group, acting properly discontinuously on the 
hyperbolic plane.
Then $\Gamma_P$ acts with finite stabilizers, and the quotient of the hyperbolic plane by $\Gamma_P$ is 
a hyperbolic orbifold $X$. We must prove that there is a compact subset of $X$ such that its intersection
with the image of any line
$L_i=S_i\cap  {\Bbb P}^+ P$ is non-empty.
But any arithmetic lattice has a finite index subgroup 
which is torsion-free (for instance, the congruence subgroup formed by integer matrices which are identity
modulo $N$ for $N$ big enough).
Therefore, our orbifold $X$ has a finite covering $\tilde{X}$ which is a hyperbolic Riemann surface,
and it suffices to prove that there is a compact $\tilde{K}\subset \tilde{X}$ such that $\pi(L_i)$ 
intersects $\tilde{K}$ for any $i$, where $\pi: {\Bbb P}^+ P\arrow \tilde{X}$ denotes the projection
(quotient by a finite index subgroup  $\tilde{\Gamma_P}\subset \Gamma_P$). 

%and \ref{_intersect_compact_Theorem_} is reduced 
%to the following lemma.

%\hfill

%%%%%%%%%%%%%%%%%%%%%%%%%%%%%%%%%%%%%%%%%%%%%%%%%%%%%
%\lemma\label{_intersect_compact_dim_2_Theorem_}
%Let $\{S_i\}$ be a set of rational lines in 
%a 2-dimensional hyperbolic space ${\Bbb P}^+ P$,
%and $\Gamma\subset SO(P, \Z)$ an arithmetic lattice.
%Then there exists a compact subset
%$K\subset {\Bbb P}^+ P$ such that $\Gamma_P\cdot K$ intersects
%all the lines $S_i$.

%\hfill

%{\bf Proof:} Denote by $\pi$ the projection
%${\Bbb P}^+ P\arrow {\Bbb P}^+ P/\Gamma=X$.
%We consider $X$ as a complete orbifold of constant negative curvature.
Let $\Gamma_{L_i}$ be the stabilizer of $L_i$ in $\Gamma_P$.
Since $\Gamma_{L_i}$ has finite index in $SO(L_i,\Z)$,
the images $\pi(L_i)$ have finite length in $\tilde{X}$.
On the other hand, $\pi(L_i)$ are isometric images
of $L_i$. Therefore, $\pi(L_i)$ are compact; in other words, these are closed geodesics on $\tilde{X}$. We have reduced
\ref{_intersect_compact_Theorem_} to the following well-known lemma.

\hfill

%%%%%%%%%%%%%%%%%%%%%%%%%%%%%%%%%%%%%%%%%%
\lemma\label{_geodesics_Lemma_}
Let $S$ be a complete hyperbolic Riemann
surface (of constant negative curvature and finite volume). Then there
exists a compact subset $K\subset X$ intersecting
each closed geodesic $l\subset S$.

\hfill

{\bf Proof:} To obtain $K$, it 
suffices to remove from $S$ a neighbourhood of each cusp: indeed,
there are no closed geodesics around cusps.
This is an elementary exercise, apparently well known;
see e.g. \cite[Theorem 1.2]{_McShane_Rivin_} which is in the same spirit. 
For the convenience of the
reader, we sketch an argument here.

Let ${\Bbb H}=\{x\in \C \ \ |\ \ \Im(x)>0\}$ 
be a hyperbolic half-plane, equipped with a Poincare metric, $t>0$
a real number, and ${\Bbb H}_t=\{x\in \C \ \ |\ \ t> \Im(x)>0\}$ 
a strip consisting of all $x\in \C$ with $0<\Im(x) <t$.
In a neighbourhood of a cusp point, $S$ is isometric to a quotient 
${\Bbb H}_t/\Z$, where the action of $\Z$ 
is generated by the parallel transport
$\gamma_r(x)= x+ r$, where $r\in \R$ is a fixed number.
A geodesic is a half-circle perpendicular to the line $\Im x=0$;
closed geodesic in ${\Bbb H}_t/\Z$ is a half-circle
which is mapped to itself by a power of $\gamma_r$. Such half-circles
clearly do not exist (see the picture).

\begin{center}\includegraphics[width=0.60\textwidth]{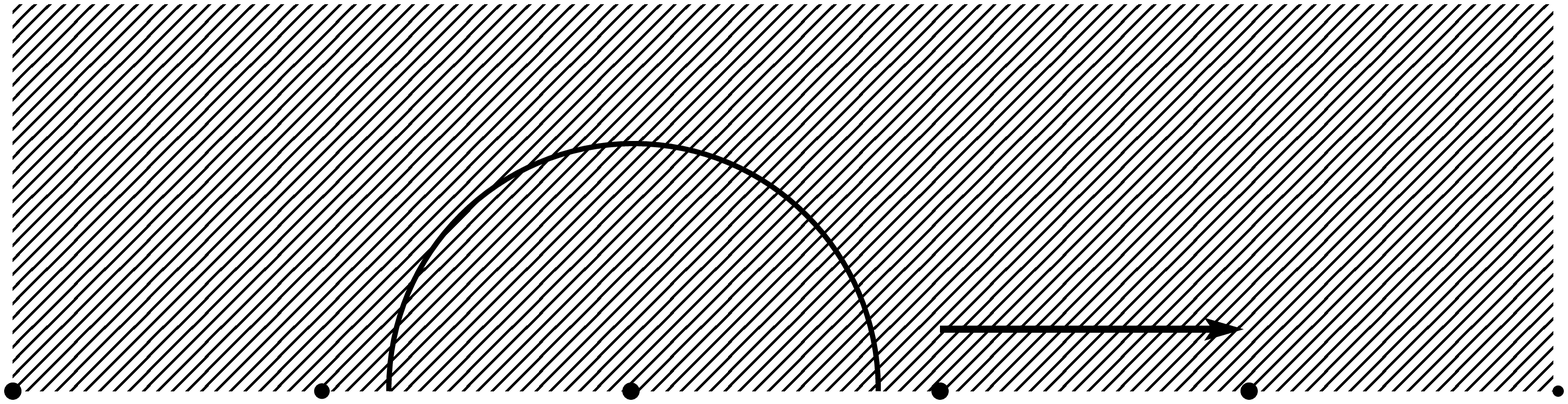}\\
{\bf \scriptsize A neighbourhood of a cusp point in dimension 2} 
\end{center}
\endproof

\hfill

\remark The result of this subsection shall be used in the next one to justify that a certain sequence of  
ergodic measures does not have a subsequence going to infinity. Since all the measures in question come from
orbits of the same subgroup, this is also a consequence of \cite{EMS}, Corollary 1.10. We prefer nevertheless
to keep our simple observations on hyperbolic geometry which might have some independent interest.

%%%%%%%%%%%%%%%%%%%%%%%%%%%%%%%%%%%%%%%%%%%%%%%%%%%%%%%%%%%%
\subsection{Measures and rational hyperplanes in the hyperbolic space}
%%%%%%%%%%%%%%%%%%%%%%%%%%%%%%%%%%%%%%%%%%%%%%%%%%%%%%%%%%%%

%%%%%%%%%%%%%%%%%%%%%%%%%%%%%%%%%%%%%%%%%%%%%%%%%%%%%%%%%%%%

The hyperbolic space, that is, the projectivization of the positive cone in a real vector space
with a quadratic form of signature $(1,n)$, is a homogeneous space in an obvious way. Indeed
it is an orbit of any positive line by the connected component of the unity of $SO(1,n)$, and the
stabilizer is isomorphic to $SO(n)$. If $z$ is a negative vector, then $z^{\perp}$ is a hyperplane
which intersects the positive cone; as in the previous paragraph, by
a hyperplane in the hyperbolic space we shall mean the projectivization of this intersection.

\hfill

%%%%%%%%%%%%%%%%%%%%%%%%%%%%%%%%%%%%%%%%%%%%%%%%%%%%% 
\theorem\label{_density_of_measures_SO(1,n)_Theorem_}
Let $G$ be the connected component of the unity $SO^+(1,n)$ of $SO(1,n)$, where $n
\geq 3$, $H:=SO(n)$, and $\Gamma\subset G_\Z$ 
a discrete subgroup of finite index 
(and therefore of finite covolume, \ref{_Borel_HCh_Remark_}).
Consider the hyperbolic space  $\H=H\backslash G=SO(n)\backslash SO^+(1,n)$.
Let $Y:= \bigcup S_i$ be a $\Gamma$-invariant union of 
rational hyperplanes. Then either $\Gamma$ acts on $\{S_i\}$
with finitely many orbits, or $Y$ is dense in $\H$.

\hfill

{\bf Proof:} 
Let $V=\R^{1,n}$ be a real vector space of signature $(1,n)$,
$G=SO^+(V)$, and $\H=H\backslash G$, where $H\subset G$ is the stabilizer
of an oriented positive hyperplane.
We may identify $\H$ with the space of positive vectors $x\in V$,
$(x,x)=1$. In order to apply ergodic theory,
we replace $\H$ by the incidence variety $X$ of pairs
$(\H_W\subset \H, x\in \H_W)$, where $\H_W$ is an oriented hyperplane
in the hyperbolic space and $x\in \H_W$. Clearly, a point of $X$ is 
uniquely determined by a pair of orthogonal
vectors $x,y\in V$, where $x$ is positive, $(x,x)=1$, and
$(y,y)=-1$. Therefore, $X=H_0 \backslash G$, where $H_0=SO(n-1)$. The important point is 
that $X$ is a quotient of $G$ by a compact group (and so is $\H$). Moreover $X$ is fibered over $\H$
in spheres of dimension $n-1$.

We can lift our hyperplanes $S_i$ to $X$ in the tautological way. To make the picture transparent,
we first treat the case $n=2$, where $\H$ is a hyperplane (since in the theorem we have $n\geq 3$,
this is just to describe the lifting and see a certain well-known analogy). Here $X=SO^+(1,2)$
is the unit tangent bundle over $\H$, and a point $x\in S_i$ lifts as $(x, z)$ where $z$ is the unit
tangent vector to $S_i$ in the direction given by the orientation. If we lift all possible 
hyperplanes to $X$ in this way, we obtain a foliation known as {\bf the geodesic flow}: our liftings never 
intersect and 
are tangent to
an invariant vector field on $X=SO^+(1,2)$. Therefore all the lifting are orbits of a Lie 
subgroup $H_1\subset G$ (this one-parameter subgroup, isomorphic to $SO^+(1,1)$, can be identified
with the group of diagonal two-by-two matrices with $e^t$ and $e^{-t}$ on the diagonal under 
an isomorphism between $SO^+(1,2)$ and $PSL(2,\R)$, see the first 
chapter of Morris' book \cite{_Morris:Ratner_}).

For $n\geq 3$, we first tautologically lift the hyperplanes to $X$ and then take preimages under the
projection from $G$ to $X$. Again, we obtain a translation-invariant foliation on $G$, which means
that the liftings and their preimages are orbits of a subgroup $H_1\subset G$ (containing $H_0$).
This subgroup is isomorphic to $SO^+(1, n-1)$, that is, generated by unipotents, so that ergodic
theory applies. 

Let us denote by $R_i$ the preimage in $G$ of the lifting of $S_i$ to $X$. Each $R_i$ is an
orbit of $H_1$. By \ref{_intersect_compact_Theorem_}, 
there is a compact
set $C$ such that the $\Gamma$-orbit 
of any $S_i$ intersects $C$. Since the projection from $G$ to $\H$ is proper, the same is true for the set
of $R_i$. Suppose that $\Gamma$ acts on the set of $S_i$ (and thus $R_i$) with infinitely many orbits.
Consider the homogeneous space $G/\Gamma$. Each $\Gamma$-orbit on the set of $R_i$ corresponds to an
algebraic probability measure $\mu_i$ on $G/\Gamma$ (note that since the hyperplanes $S_i$ are rational,
the quotient of each $H_1$-orbit $R_i$ over its stabilizer in $\Gamma$ has finite Haar volume by Borel and
Harish-Chandra theorem). The support of $\mu_i$ is the image of $R_i$ in $G/\Gamma$. Since the union
of $R_i$ is $\Gamma$-invariant, to prove 
\ref{_density_of_measures_SO(1,n)_Theorem_},
it suffices to show that the union of $\Supp(\mu_i)$ is dense in $G/\Gamma$:
this will imply the density of $R_i$ in $G$ and therefore the density of $S_i$ in $\H$. 

By \ref{_compactness_MS_Theorem_}, the sequence $\mu_i$ has a limit point which is either a 
probability measure, or is supported at infinity. But the latter option is impossible. Indeed, by 
\ref{_intersect_compact_Theorem_}
all $\Supp(\mu_i)$ intersect the same compact on $G/\Gamma$. Thus there is a (slightly larger) compact $C$
such that $\mu_i(C)>0$ for all $i$, and by Dani-Margulis theorem, for another compact $K_{\epsilon}$ and 
all $i$, $\mu_i(K)>1-\epsilon$.

Taking a suitable subsequence, we may therefore suppose that $\lim \mu_i=\mu$ where $\mu$ is an
algebraic probability measure.

%Consider the subgroup $H_1\supset H_0$ of $G$ fixing $y$.
%By definition, orbits of $H_1$ in $X$ are all 
%pairs $(x\in \H, y)$, where $x$ is 
%contained in a given hyperplane $\H_W=y^\bot$,
%and $y$ fixed. This construction
%is used to identify
%the set of hyperplanes in $\H$ with orbits of $H_1\cong SO^+(1,n-1)$
%acting on $X=G/H_0$.

%To prove \ref{_density_of_measures_SO(1,n)_Theorem_},
%consider the set $M_1$ of algebraic probability
%measures $\mu$ on $\Gamma\backslash \H$ with $\Supp(\mu)$ 
%equal to a hyperplane $\H_W\subset \H$.
%Clearly, $M_1$ is in bijective correspondence with the set
%of rational hyperplanes in $\H$ up to $\Gamma$-action.
%However, the projection $X\arrow \H$ is proper. Therefore,
%to prove that $Y$ is dense, it would suffice to show that
%the corresponding union of $H_1$-orbits is dense in $X$.

We have reduced  \ref{_density_of_measures_SO(1,n)_Theorem_}
to the following lemma.

\hfill

%%%%%%%%%%%%%%%%%%%%%%%%%%%%%%%%%%%%%%%%%%%%%%%%%%%%% 
\lemma\label{_same_subgroup_}
Let $G$ be the connected component $SO^+(1,n)$ of $SO(1,n)$, where $n
\geq 3$, and $\Gamma\subset G_\Z$ 
a discrete subgroup of finite index (and therefore of finite covolume).
Let $H_1\subset G$ be $SO^+(1,n-1)$.
Let $\mu_i$ be a sequence of algebraic probability measures on $G/\Gamma$
associated with the orbits of $H_1$. Suppose $\mu_i$ converges to an 
algebraic probability measure $\mu$. Then either $\mu_i$ are finitely many, or 
$\Supp(\mu)$ is $G/\Gamma$, so that $\Supp(\mu_i)$ are dense in $G/\Gamma$.
%Consider the space $X=G/H_0=SO(1,n)/SO(n-1)$.
%Let $Y:= \bigcup S_i$ be a $\Gamma$-invariant union of 
%rational $H_1$-orbits. Then either $\Gamma$ acts on $\{S_i\}$
%with finitely many orbits, or $Y$ is dense in $X$.

\hfill

%{\bf Proof. Step 1:} The space of probability measures 
%on the one-point compactification of
%$G/\Gamma$ is compact by \ref{_compactness_MS_Theorem_}. 
%Let ${\cal P}$ be the space of $H_0$-invariant
%probability measures 
%on the one-point compactification of
%$\Gamma\backslash G$, or, equivalently, 
%the space of probability measures 
%on the one-point compactification of
%$\Gamma\backslash X$. Clearly, ${\cal P}$ 
%is also compact.

%Each rational $H_1$-orbit $S_i$ gives a point in ${\cal P}$.
%Indeed, the quotient of $S_i$ over its 
%stabilizer in $\Gamma$ has finite Haar volume by
%Borel Harish-Chandra theorem (\ref{_Borel_HCh_Remark_}).

%Suppose that $\Gamma$ acts on the
%set of $S_i$ with infinitely
%many orbits, giving infinitely many distinct
%measures $\mu_i\in {\cal P}$. Let $\mu$ be any limit point of the
%sequence $\{\mu_i\}$. By 
%\ref{_intersect_compact_Theorem_}, supports $R_i$ of $\mu_i$
%all intersect a compact subset. Therefore, 
%by Dani-Margulis theorem (\ref{_Dani_Margulis_Theorem_}), 
%$\mu$ is not supported at infinity and so induces an
%algebraic probability measure on $\Gamma\backslash X$.

%{\bf Step 2:} Taking a subsequence, we may assume
%that $\lim\mu_i=\mu$. 

{\it Proof:} By \ref{_MS_main_Theorem_} (ii), the
support of $\mu$ contains a right translate by $g_i\to e$ 
of the support of infinitely many of $\mu_i$. 
Moreover, $\mu$ is an algebraic measure associated with an 
orbit of a closed subgroup $F\subset G$. 
But there are no closed intermediate connected subgroups between 
$G=SO^+(1,n)$ and $H_1$, which stabilizes 
a hyperplane. Therefore, $F$ is either equal to $G$, or
is the stabilizer $H$ of a hyperplane $\H_W$.

In the first case, the support of $\mu=\lim\mu_i$ is
$G/\Gamma$ and thus $\Supp(\mu_i)$ are dense in
$G/\Gamma$. 

In the second case, for $i\gg 0$, $\Supp(\mu)=g_i\Supp(\mu_i)$,
that is, $Hx=g_iH_1x_i$ where $g_ix_i=x$ (where $x$, $x_i$, $g_i$ are as in \ref{_MS_main_Theorem_}).
That is, $Hx=H_1^{g_i}x$ and therefore $H_1^{g_i}=H=H_1$.
Since $H_1$ has finite index in its normalizer, this means that
there are only finitely many 
$\mu_i$.
%right translations $R_{g_i}$ by $g_i$
%map the $H_1$-orbit $S_i=x_i\cdot H_1$ to the right orbit 
%$\Supp(\mu)=x\cdot H_1$:
%\[ R_{g_i}(x_i\cdot H_1)= x\cdot H_1= x_i g_i H_1^{g_i}.
%\]
%This gives $x=x_i g_i$ and $H_1^{g_i}=H_1$,
%hence all $S_i$ are equal for $i\gg 0$.
\endproof

%%%%%%%%%%%%%%%%%%%%%%%%%%%%%%%%%%%%%%%%%%%%%%%%%%%%%%%%%%%%

\section{The proof of Morrison-Kawamata cone conjecture}

%%%%%%%%%%%%%%%%%%%%%%%%%%%%%%%%%%%%%%%%%%%%%%%%%%%%%%%%%%%%

%%%%%%%%%%%%%%%%%%%%%%%%%%%%%%%%%%%%%%%%%%%%%%%%%%%%%%%%%%%%
\subsection{Morrison-Kawamata conjecture for the K\"ahler cone}
%%%%%%%%%%%%%%%%%%%%%%%%%%%%%%%%%%%%%%%%%%%%%%%%%%%%%%%%%%%%

The following theorem is an immediate consequence of \ref{_density_of_measures_SO(1,n)_Theorem_}.

\hfill

\theorem\label{finitely_many_orbits} Let $M$ be a projective simple hyperk\"ahler manifold which has
Picard number at least 4. Then the Hodge monodromy group acts with 
finitely many orbits on the set of MBM classes of type $(1,1)$.

\hfill

{\bf Proof:} This is the same as to say that the Hodge monodromy group acts with finitely many orbits on the set
of their orthogonal hyperplanes, which by \cite{_AV:MBM_} are exactly the hyperplanes supporting the faces of the
K\"ahler chambers.

Since the Hodge monodromy group is of finite index in the orthogonal group of the Picard lattice, which is
of signature $(+,-,\dots, -)$, one can apply
\ref{_density_of_measures_SO(1,n)_Theorem_} to the Picard lattice, with $\Gamma$
the Hodge monodromy group. One concludes that if the number of
$\Gamma$-orbits is infinite, then the hyperplanes orthogonal to MBM classes should be dense in the positive
cone. This is clearly absurd, as they should bound 
the K\"ahler cone (Subsection \ref{_Kahler_cone_MBM_intro_Subsection_}), 
so the number of $\Gamma$-orbits 
is finite. \endproof

\hfill

\corollary On an $M$ as above, the
primitive MBM classes of type $(1,1)$ have bounded 
Beauville-Bogomolov square.

\hfill

{\bf Proof:} Indeed, the monodromy acts by isometries.
\endproof

\hfill

\theorem\label{boundedness-nonproj} Let $M$ be a
simple hyperk\"ahler manifold such that $b_2(M)\geq 6$.
Then the primitive MBM classes of type $(1,1)$ have bounded 
Beauville-Bogomolov square.

\hfill

{\bf Proof:} If $M$ is not projective or the Picard number of $M$ is less than four, apply \ref{defoproj} 
to get a projective deformation $M'$ with Picard number at least four such that
all MBM classes of type $(1,1)$ on $M$ remain of type $(1,1)$ on $M'$. Then use the deformation invariance of
MBM property proved in \cite{_AV:MBM_} to conclude that these MBM classes remain MBM on $M'$ and therefore the primitive
ones must have bounded square by the preceding theorem.

\hfill

The Morrison-Kawamata conjecture for the K\"ahler cone
now follows in the same way as in \cite{_AV:MBM_}.

\hfill

\theorem\label{cone} Let $M$ be a simple hyperk\"ahler manifold with $b_2(M)\geq 6$. Then the automorphism group of $M$ acts with
finitely many orbits on the set of faces of its K\"ahler cone.

\hfill

{\bf Proof:} The argument is the same as in \cite{_AV:MBM_} where the theorem has been obtained under the
boundedness assumption on squares of primitive MBM classes, which we have just proved: see \cite[Theorem 6.6]{_AV:MBM_} there for an outline of the
argument and \cite[Theorem 3.14, 3.29]{_AV:MBM_} for technicalities. \endproof

%%%%%%%%%%%%%%%%%%%%%%%%%%%%%%%%%%%%%%%%%%%%%%%%%%%%%%%%%%%%%%%%%%%%
\subsection{Morrison-Kawamata conjecture for the ample cone}
\label{_ample_Subsection_}
%%%%%%%%%%%%%%%%%%%%%%%%%%%%%%%%%%%%%%%%%%%%%%%%%%%%%%%%%%%%%%%%%%%%

Recall from e.g. \cite{MY} that the classical Morrison-Kawamata cone conjecture is formulated in
the projective case and treats the ample cone rather than the K\"ahler cone. It also states something 
{\em a priori} 
stronger 
than the finiteness of the number of orbits of the action of automorphism group on the set of faces of the
cone, namely the existence of a finite polyhedral fundamental domain.

\hfill

\conjecture (Morrison-Kawamata cone conjecture for the ample cone)\\
\label{mor-kawbis} 
The automorphism group $\Aut(M)$ has a finite polyhedral
fundamental domain on the ample cone.

\hfill

We shall see in this subsection that this in fact follows from our version of the cone conjecture,
and therefore is true for all simple hyperk\"ahler manifolds with $b_2\neq 5$.

The ample cone $\Amp(M)$ is the convex hull of $\Kah(M)\cap H^{1,1}_{\Q}(X)$ in the space 
$H^{1,1}_{\Q}(X)\otimes \R=
\NS(X)\otimes \R$, so that 
{\em a priori} new faces could arise from the ``circular part'' of the boundary of the K\"ahler cone. In our case,
this is not a problem, since this part is a piece of the quadric over the rationals 
defining $\Pos(M)$, and when it has a single rational point, it has a dense set of them. Thus the
Hodge monodromy group acts with finitely many orbits on the set of faces of the ample cone.

Denote by ${\cal C}(M)$ the intersection of $\Pos (M)$ with $\NS(X)\otimes \R$. The Hodge monodromy
group $\Gamma$ acts on $\P{\cal C}(M)$ with finite stabilizers (since the stabilizer of a point $x$ in
$\P{\cal C}(M)$ must also stabilize the orthogonal hyperplane to the line corresponding to $x$,
and our form is negative definite on such a hyperplane). By its arithmeticity,
replacing if necessary the group $\Gamma$ by a finite index subgroup, we may assume there are
no stabilizers at all.
Indeed, an arithmetic lattice has a finite index
torsion-free subgroup, which can be obtained by taking
a congruence subgroup formed by integer matrices which are identity
modulo $N$ for $N$ big enough.
Consider the quotient $S:={\cal C}(M)/\Gamma$. Since
$\Gamma$ is arithmetic, Borel and Harish-Chandra
theorem implies that $S$ is a complete hyperbolic
manifold of finite volume. The image of $\Amp(M)$ in $S$ is a hyperbolic 
manifold $T$ with finite (that is, consisting of finitely many geodesic pieces)
boundary, by \ref{cone}. It is known 
(see \cite[Proposition 4.7 and 5.6]{bowditch} 
or \cite[Theorem 2.6]{kapovich})  that such
manifolds are {\it geometrically finite}, that is, they admit a finite cell decomposition with finite
piecewise geodesic boundary (in fact one even has a decomposition with a single cell of maximal dimension,
the {\it Dirichlet-Voronoi
decomposition}). Thus $T$ becomes a union of finitely many cells with finite
piecewise geodesic boundary. Taking the union of suitable liftings of these to $\Amp(M)$,
we obtain a finite polyhedron within the closure of $\Amp(M)$ which is a fundamental domain for the
subgroup of $\Gamma$ preserving $\Amp(M)$, that is, of the automorphism group of $M$. We thus have
proved the following

\hfill

\theorem Let $M$ be a projective simple hyperk\"ahler manifold with $b_2\neq 5$. The automorphism
group has a finite polyhedral fundamental domain on the ample cone of $M$.

\hfill

{\bf Acknowledgements:}
We are
grateful to Eyal Markman for many interesting discussions
and an inspiration.  Many thanks to Alex Eskin,
Maxim Kontsevich, Anton Zorich, Vladimir Fock and S\'ebastien Gouezel for elucidating the various
points of measure theory and hyperbolic geometry. The scheme of the proof 
of \ref{_density_of_measures_SO(1,n)_Theorem_} is due to
Alex Eskin; we are much indebted to Alex for his
kindness and patience in explaining it. Much gratitude to
Mihai Paun and Sasha Anan'in for inspiring discussions. Thanks to Martin Moeller for indicating the reference \cite{EMS}, and to Andrey
Konovalov for solving a tricky problem of 
hyperbolic plane geometry which disproved 
one of our conjectures.

{%\scriptsize

\noindent {\sc Ekaterina Amerik\\
{\sc Laboratory of Algebraic Geometry,\\
National Research University HSE,\\
Department of Mathematics, 7 Vavilova Str. Moscow, Russia,}\\
\tt  Ekaterina.Amerik@math.u-psud.fr}, also: \\
{\sc Universit\'e Paris-11,\\
Laboratoire de Math\'ematiques,\\
Campus d'Orsay, B\^atiment 425, 91405 Orsay, France}

\hfill

\noindent {\sc Misha Verbitsky\\
{\sc Laboratory of Algebraic Geometry,\\
National Research University HSE,\\
Department of Mathematics, 7 Vavilova Str. Moscow, Russia,}\\
\tt  verbit@mccme.ru}, also: \\
{\sc Kavli IPMU (WPI), the University of Tokyo}

}

\end{document}